\documentclass[11pt,a4paper]{article}
\pdfoutput=1

\usepackage{amsfonts,amsmath,amssymb,graphicx}
\usepackage[a4paper,margin=1.2in]{geometry}

\newtheorem{theorem}{\bf Theorem}[section]
\newtheorem{corollary}[theorem]{\bf Corollary}
\newtheorem{lemma}[theorem]{\bf Lemma}
\newtheorem{proposition}[theorem]{\bf Proposition}

\newcommand{\proof}{\noindent{\bf Proof.\ }}
\newcommand{\qed}{\hfill $\Box$ \bigskip}

\newcommand{\name}{X}
\newcommand{\idim}{{\rm idim}}
\newcommand{\fdim}{{\rm fdim}}
\newcommand{\ldim}{{\rm ldim}}
\def\cp{\,\square\,}
\newcommand{\W}{{\mathbb{W}}}
\newcommand{\simplex}{\kappa}

\begin{document}
\title{The Fibonacci dimension of a graph}

\author{Sergio Cabello\thanks{%
	Faculty of Mathematics and Physics, University of Ljubljana,
	Jadranska 19, 1000 Ljubljana, Slovenia; 
	Institute of Mathematics, Physics and Mechanics, Jadranska 19, 1000 Ljubljana, Slovenia.
	E-mail: {\tt sergio.cabello@fmf.uni-lj.si}.
	}
\and David Eppstein\thanks{%
	Computer Science Department,
	University of California, Irvine, CA 92697-3425, USA.
	Email: {\tt eppstein@uci.edu}.}
\and
    Sandi Klav{\v z}ar\thanks{%
    Faculty of Mathematics and Physics, University of Ljubljana,
	Jadranska 19, 1000 Ljubljana, Slovenia;
	Faculty of Natural Sciences and Mathematics, University of Maribor,
	Koro{\v s}ka 160, 2000 Maribor, Slovenia;
	Institute of Mathematics, Physics and Mechanics, Jadranska 19, 1000 Ljubljana, Slovenia.
	E-mail: {\tt sandi.klavzar@fmf.uni-lj.si}. }
}
\date{\today}

\maketitle

\begin{abstract}
The Fibonacci dimension $\fdim(G)$ of a graph $G$ is introduced as the 
smallest integer $f$ such that $G$ admits an isometric embedding into 
$\Gamma_f$, the $f$-dimensional Fibonacci cube.
We give bounds on the Fibonacci dimension of a graph
in terms of the isometric and lattice dimension, 
provide a combinatorial characterization 
of the Fibonacci dimension using properties of an associated graph, 
and establish the Fibonacci dimension for certain families of graphs.
From the algorithmic point of view we prove that it is NP-complete 
to decide if $\fdim(G)$ equals to the isometric dimension of $G$, and that 
it is also NP-hard to approximate $\fdim(G)$ within $(741/740)-\varepsilon$.
We also give a $(3/2)$-approximation algorithm for $\fdim(G)$ in the general 
case and a $(1+\varepsilon)$-approximation algorithm for simplex graphs. 
\end{abstract}

%

\section{Introduction}

Hypercubes play a prominent role in metric graph theory as well as 
in several other areas such as parallel computing and coding theory. One of their 
central features is the ability to compute distances very efficiently
because the distance between two vertices is simply the number of coordinates
in which they differ; the same ability to compute distances may be transferred to any isometric subgraph of a hypercube. In this way partial cubes
appear, a class of graphs intensively studied so far;
see the books~\cite{dela-97,epfaov-08,imkl-00}, 
the recent papers~\cite{begr-08,pize-08,po-07,po-08}, 
the recent (semi-)survey~\cite{ov-08}, and references therein. In particular we point 
out a recent fast recognition algorithm~\cite{ep-08} 
and
improvements in classification of cubic partial cubes~\cite{ep-06,klsh-07}. 

The isometric dimension of a graph $G$ is the smallest (and at the same time the 
largest) integer $d$ 
such that $G$ isometrically and irredundantly embeds into the $d$-dimensional
cube. Clearly, the isometric dimension of $G$ is finite if and only if
$G$ is a partial cube. This graph dimension is well-understood; for 
instance, it is equal to the number of steps in Chepoi's expansion 
procedure~\cite{ch-88} and to the number of $\Theta$-equivalence 
classes~\cite{dj-73,wi-84} of a given graph. 
Two related graph dimensions need to be mentioned here since they are both 
defined on the basis of isometric embeddability into graph products. 
The lattice dimension of a graph is the smallest $d$ such that the
graph embeds isometrically into $\mathbb{Z}^d$. Graphs with finite
lattice dimension are precisely partial cubes and the dimension
can be determined in polynomial time~\cite{ep-05}. Another dimension
is the strong isometric dimension---the smallest integer $d$ 
such that a graph isometrically embeds into the strong product of 
$d$ paths~\cite{fino-00,frje-07}. In this case every 
graph has finite dimension, but this universality has a price: it is
very difficult to compute the strong isometric dimension.  

Fibonacci cubes were first introduced by Hsu et al. in 1993~\cite{hsu-93,hsu-93b}, although closely related structures had been studied previously~\cite{Bec-FQ-90,Gan-DM-82,HofHof-FQ-85}. Different structural properties 
of this class of graphs were investigated~\cite{deto-02,klzi-05,muza-02}.  
In~\cite{brkl-06} it was shown that Fibonacci cubes are $\Theta$-graceful
while in~\cite{tave-07} an efficient recognition algorithm is 
presented. The original motivation for introducing Fibonacci cubes was as an interconnection network for parallel computers; in that application, it is of interest to study the embeddability of other networks within Fibonacci cubes~\cite{ConZheSha-P7IPPS-93,GonZhe-SESST-96}.

In this paper we study this embedding question from the isometric point of view. 
We introduce the Fibonacci dimension of a graph as the 
smallest integer $f$ such that the graph admits an isometric 
embedding into the $f$-dimensional Fibonacci cube. In the next section 
we give definitions, notions, and preliminary results needed
in this paper. 
In Section~\ref{sec:combin} we a give a combinatorial characterization 
of the Fibonacci dimension using properties of an associated graph,
provide upper and lower bounds for the Fibonacci dimension in terms
of the isometric and lattice dimension, 
and discuss the Fibonacci dimension of some particular classes of graphs.
In Section~\ref{sec:algor} we show that computing the Fibonacci dimension
is an NP-complete problem, provide inapproximability results,
and give approximation algorithms.

\section{Preliminaries}

\begin{figure}[t]
\centering\includegraphics[width=3in]{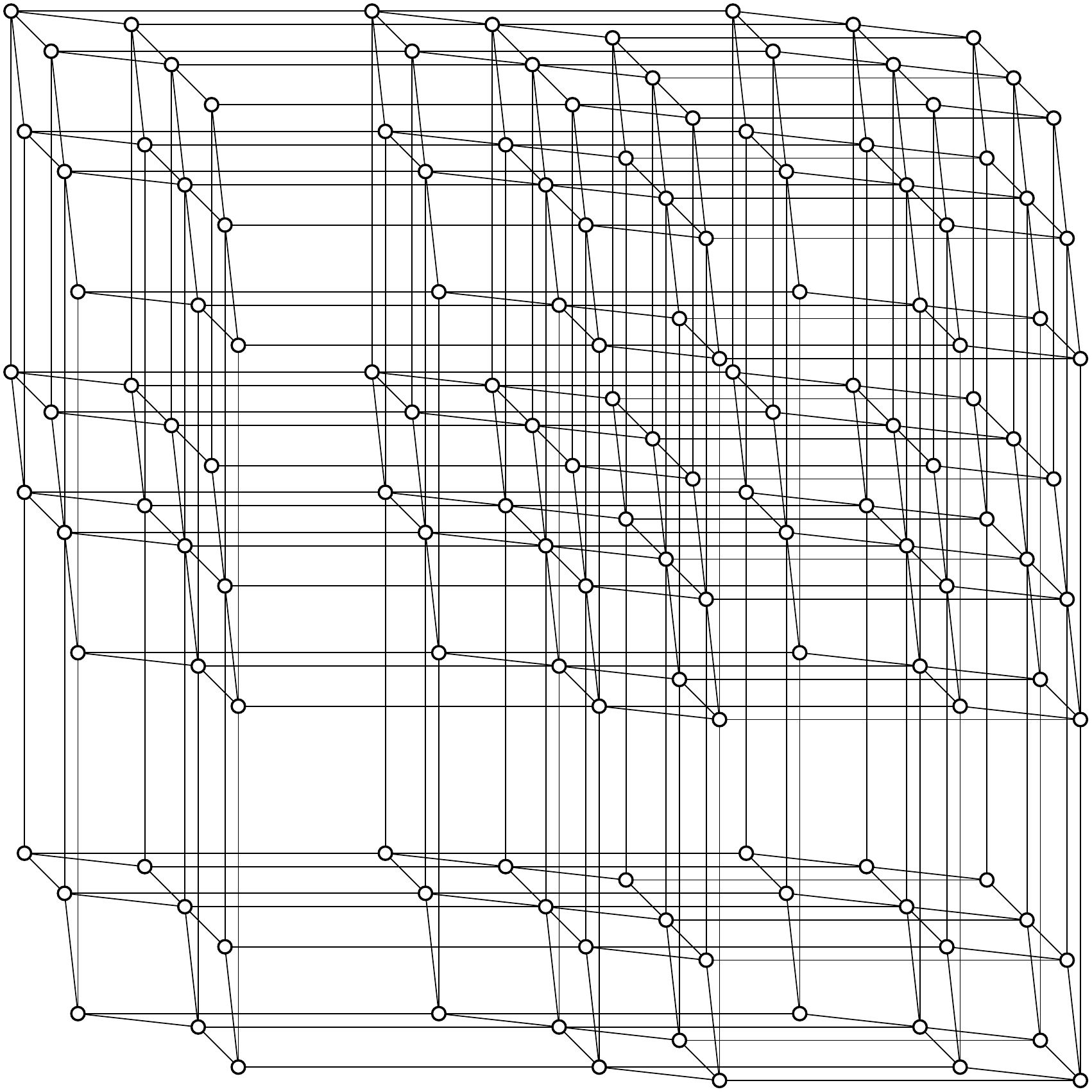}
\caption{The Fibonacci cube $\Gamma_{10}$.}
\end{figure}

We will use the notation $[n]=\{1,\ldots,n\}$.
For any string $u$ we will use $u^{(i)}$ to denote its $i$th coordinate.
Unless otherwise specified, the distance in this paper is the usual shortest-path distance for unweighted graphs. A graph $G$ is an \emph{isometric subgraph} of another graph $H$ if there is a way of placing the vertices of $G$ in one-to-one correspondence with a subset of vertices of $H$, such that the distance in $G$ equals the distance between corresponding vertices in $H$.

The vertex set of the \emph{$d$-cube} $Q_d$ 
consists of all $d$-tuples $u=u^{(1)}u^{(2)}\ldots u^{(d)}$ with 
$u^{(i)} \in \{0,1\}$. Two vertices are adjacent if the 
corresponding tuples differ in precisely one position. $Q_d$ is also 
called a \emph{hypercube of dimension} $d$.
Isometric subgraphs of hypercubes are \emph{partial cubes}. 

A \emph{Fibonacci string} of length $d$ is a binary string
$u^{(1)}u^{(2)}\ldots u^{(d)}$ with $u^{(i)}\cdot u^{(i+1)}=0$ for $i\in [d-1]$. In
other words, a Fibonacci string is a binary string without two 
consecutive ones. The set of Fibonacci strings of length $d$ can be decomposed into two subsets, strings starting with $0$ followed by a Fibonacci string of length $d-1$, and strings starting with $10$ followed by a Fibonacci string of length $d-2$. For this reason the number of distinct Fibonacci strings of length $d$ satisfies the Fibonacci recurrence and equals a Fibonacci number. The \emph{Fibonacci cube} 
$\Gamma_d$, $d\geq 1$, is the subgraph of $Q_d$ induced by 
the Fibonacci strings of length $d$.
The Fibonacci cube may alternatively be defined as the graph of the distributive lattice of order-ideals of a fence poset~\cite{Bec-FQ-90,Gan-DM-82,HofHof-FQ-85} or as the simplex graph of the complement graph of a path graph. Since graphs of distributive lattices and simplex graphs are both instances of median graphs~\cite{Banvan-PotLMS-89,BirKis-BotAMS-47}, we have:

\begin{theorem}[\cite{kl-05}]
\label{thm:old-median}
	Fibonacci cubes are median graphs. In particular, 
	Fibonacci cubes are partial cubes and $\Gamma_d$ isometrically 
	embeds into $Q_d$.  
\end{theorem}

We will use the lattice $\mathbb{Z}^d$ equipped with the $L_1$-distance.
Therefore, the distance between any two elements $(x_1,\dots,x_d),(y_1,\dots,y_d)\in \mathbb{Z}^d$
is given by $\sum_i |x_i-y_i|$. It will be convenient to visualize $\mathbb{Z}^d$
as an infinite graph whose vertex set are elements of $\mathbb{Z}^d$ and where two
vertices are adjacent when they are at distance one; with this visualization, $L_1$-distance coincides with the shortest path distance in the graph.

Let $G$ be a connected graph. The \emph{isometric dimension}, $\idim(G)$, 
is the smallest integer $k$ such that $G$ admits an isometric embedding 
into $Q_k$. If there is no such $k$ we set $\idim(G) = \infty$.
Be definition, $\idim(G) < \infty$ if and only if $G$ is a 
partial cube. 
The \emph{lattice dimension}, $\ldim(G)$, 
is the smallest integer $\ell$ such that $G$ admits an isometric embedding 
into $\mathbb{Z}^\ell$.
We similarly define the \emph{Fibonacci dimension}, $\fdim(G)$, 
as the smallest integer 
$f$ such that $G$ admits an isometric embedding into $\Gamma_f$, and 
set $\fdim(G) = \infty$ if there is no such $f$.

Let $\beta: V(G)\rightarrow V(Q_k)$ be an isometric embedding. 
We will denote the $i$th coordinate of $\beta$ with $\beta^{(i)}$. 
The embedding $\beta$ is called \emph{irredundant}
if $\beta^{(i)}(V(G))=\{ 0,1 \}$ for each $i\in [k]$.
If an embedding is not irredundant, we may find an embedding onto a lower-dimensional hypercube by omitting the redundant coordinates.
An isometric embedding
$\beta: \, G\rightarrow Q_k$ is irredundant if and only if 
$k=\idim(G)$~\cite{wi-84}.

\begin{figure}[t]
\centering\includegraphics[scale=0.5]{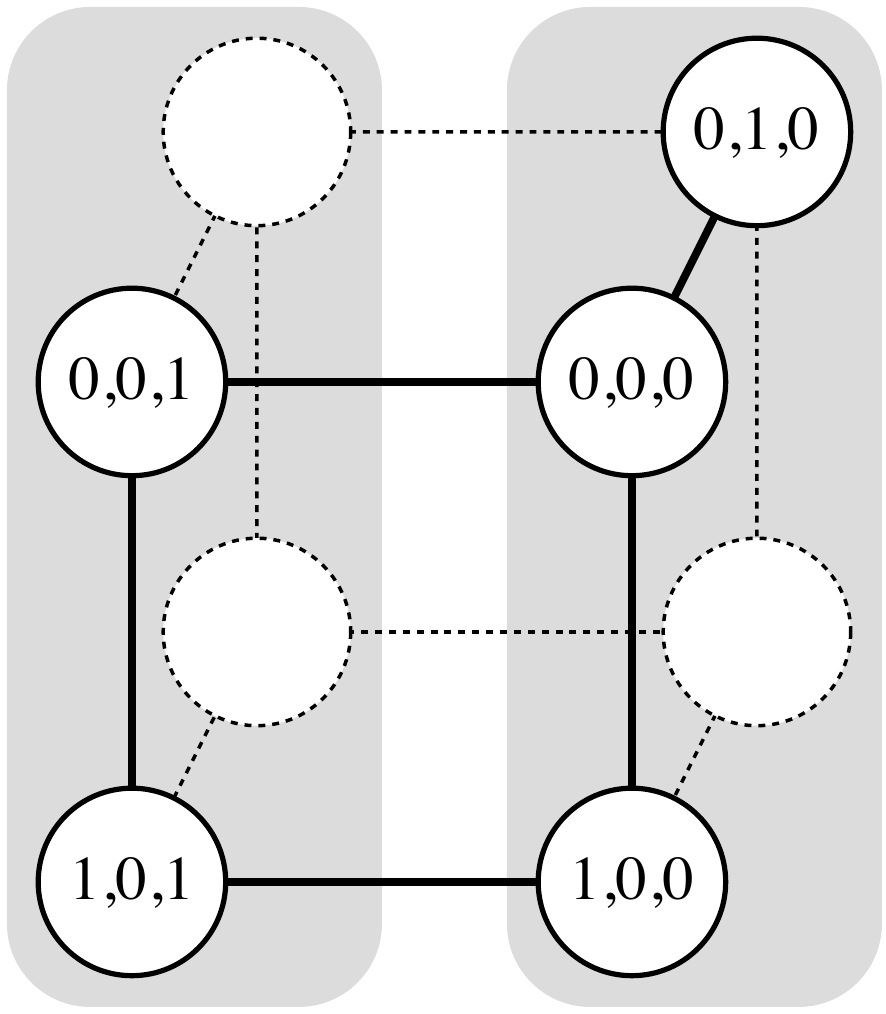}\qquad\qquad
\raise0.5in\hbox{\includegraphics[scale=0.6]{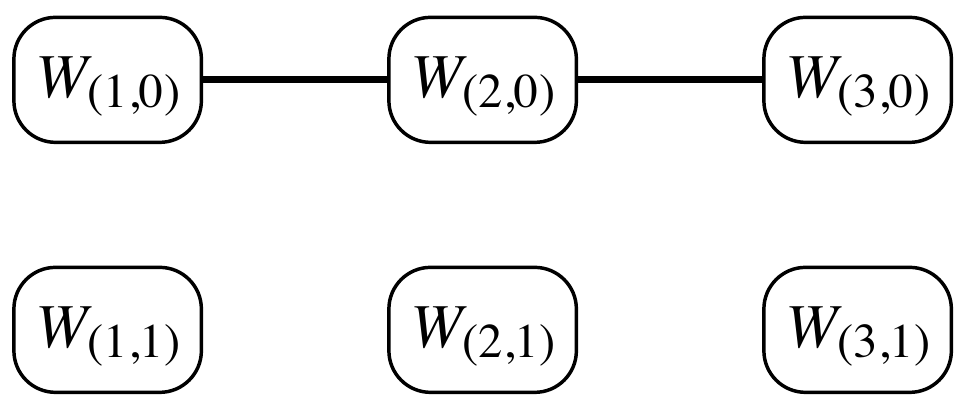}}
\caption{Left: an isometric embedding of $\Gamma_3$ into $Q_3$, with the complementary semicubes $W_{(3,1)}$ and $W_{(3,0)}$ shown as the shaded regions of the drawing. Right: the semicube graph of the embedding, consisting of a three-vertex path and three isolated vertices.}
\label{fig:Q3G3-scg}
\end{figure}

Let $G$ be a partial cube with $\idim(G)=k$ and assume that
we are given an isometric embedding $\beta$ of $G$ into $Q_k$.
Each pair $(i,\chi)\in [k]\times \{0,1\}$ defines the \emph{semicube}
$W_{(i,\chi)}=\{ u\in V(G)\mid \beta^{(i)}(u)=\chi \}$.
For any $i\in [k]$, we refer to $W_{(i,0)}, W_{(i,1)}$ as 
a \emph{complementary pair of semicubes}.
This definition and notation seems to depend on the embedding $\beta$. 
However, any irredundant isometric embedding $\beta'$ describes the same 
family of semicubes and pairs of complementary semicubes, possibly indexed
in a different way.

For a partial cube $G$ and a complementary pair of semicubes
$W_{(i,0)}, W_{(i,1)}$, the set of edges with one endvertex in 
$W_{(i,0)}$ and the other in $W_{(i,1)}$ constitute a \emph{$\Theta$-class} 
of $G$. The $\Theta$-classes of $G$ form a partition of $E(G)$. 

To determine the lattice dimension of a graph $G$, Eppstein~\cite{ep-05} 
introduced the semicube graph Sc$(G)$ of a partial cube $G$ as 
the graph with all the semicubes as nodes, semicubes 
$W_{(i,\chi)}$ and $W_{(i',\chi')}$ being adjacent 
if $W_{(i,\chi)} \cup W_{(i',\chi')} = V(G)$ and 
$W_{(i,\chi)} \cap W_{(i',\chi')} \not= \emptyset$.
One can then show that the lattice dimension of $G$ is equal to 
$\idim(G) - |M|$, where $M$ is a maximum matching of $\mbox{Sc}(G)$.
See also~\cite{klko-09} for further work on semicube graphs.

\begin{figure}
	\centering
	\includegraphics[width=\textwidth]{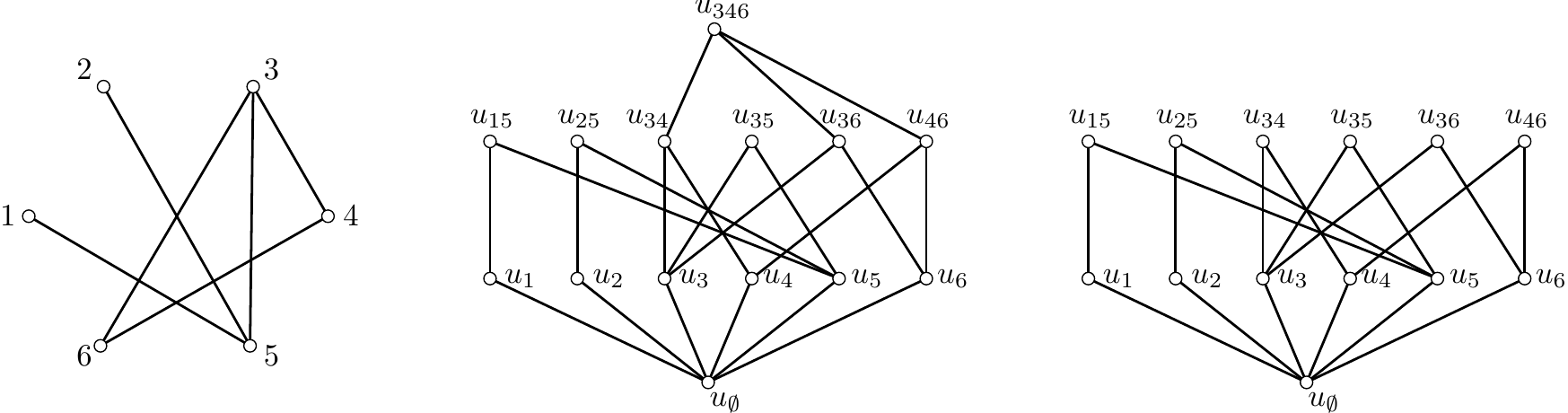}
	\caption{A graph $G$ (left) with its corresponding simplex graph $\simplex(G)$ (center) 
			and 2-simplex graph $\simplex_2(G)$.}
	\label{fig:HG}
\end{figure}

For any graph $G$, its \emph{simplex graph} $\simplex(G)$ is defined as follows.
There is a vertex $u_K$ in $\simplex(G)$ for each clique $K$ of $G$;
here we regard $\emptyset$, each vertex, and each edge of $G$ as a clique.
There is an edge between vertices $u_K$ and $u_{K'}$ of $\simplex(G)$
whenever the cliques $K$ and $K'$ of $G$ differ by exactly one vertex.
In particular, there is an edge between $u_\emptyset$ and $u_{a}$ for each $a\in V(G)$,
and there is an edge between $u_a$ and $u_{ab}$ for each edge $ab\in E(G)$.
We will also use the \emph{2-simplex graph} $\simplex_2(G)$ of a graph $G$,
which is the subgraph of $\simplex(G)$ induced by the vertices 
$u_K$ of $\simplex(G)$ corresponding to cliques $K$ with at most 2 vertices.
An example is given in Figure~\ref{fig:HG}.
When $G$ has no triangle, then
$\simplex_2(G)=\simplex(G)$.
2-simplex graphs were used in~\cite{imklmu-99} to 
establish a close connection between the recognition 
complexity of triangle-free graphs and of median graphs. 

Finally, \emph{computing an embedding} of $G$ into $Q_d$ (or $\Gamma_d$) 
means to attach to each vertex $v$ of $G$ a tuple $\beta(v)$ that 
is a vertex of $Q_d$ such that $\beta$ provides an isometric embedding.

\section{Combinatorial aspects}
\label{sec:combin}

\subsection{The general case}
\label{sec:combin-general}

\begin{proposition}
\label{prp:pc}
	Let $G$ be a connected graph. Then $\fdim(G) < \infty$ if and only 
	if $\idim(G) < \infty$. Moreover, 
	$$\idim(G) \le \fdim(G) \le 2\,\idim(G) -1\,.$$
\end{proposition}

\proof
	Let $f=\fdim(G) < \infty$, so that $G$ isometrically embeds into 
	$\Gamma_f$. By Theorem~\ref{thm:old-median}, $\Gamma_f$ isometrically
	embeds into $Q_f$, hence $G$ isometrically embeds into $Q_f$.  The Fibonacci strings with which $\Gamma_f$ was derived may be used directly as the coordinates of an isometric embedding.
	Consequently $\idim(G) \le f = \fdim(G)$. 

	Conversely, let $k = \idim(G) < \infty$ and consider $G$ 
	isometrically embedded into $Q_k$. To each vertex $u=u^{(1)}u^{(2)}\ldots u^{(k-1)}u^{(k)}$
	of $G$ (embedded into $Q_k$) assign the vertex 
	$\widetilde{u} = u^{(1)}0u^{(2)}0\ldots u^{(k-1)}0u^{(k)}$. Clearly, 
	$\widetilde{u}^{(i)}\cdot \widetilde{u}^{(i+1)}=0$ for any $i\in [2k-2]$. 
	Therefore, we can consider $\widetilde{u}$ as a vertex of $\Gamma_{2k-1}$. 
	Let $\widetilde{G}$ be the subgraph of $\Gamma_{2k-1}$ induced
	by the vertices $\widetilde{u}$, $u\in V(G)$. 
	Since $\Gamma_{2k-1}$ is isometric in $Q_{2k-1}$ (invoking 
	Theorem~\ref{thm:old-median} again), it readily follows that 
	$\widetilde{G}$ is isometric in $\Gamma_{2k-1}$. We conclude that 
	$\fdim(G) \le 2k-1 = 2\,\idim(G)-1$. 
\qed

It is now clear that we only need to study the Fibonacci dimension of partial cubes.
Using the lattice dimension $\ldim(G)$ we will further improve in Proposition~\ref{prp:ldim1}
the upper bound on $\fdim(G)$, and provide an alternative lower bound in Proposition~\ref{prp:ldim2}.

Let $G$ be a partial cube with $\idim(G) = k$. 
In order to obtain an expression for $\fdim(G)$ in terms of $\idim(G)$ 
we construct the graph $\name(G)$ as follows. 
The nodes of $\name(G)$ are the semicubes 
$W_{(i,\chi)}$, $(i,\chi)\in [k]\times \{0,1 \}$, of $G$, 
semicubes $W_{(i,\chi)}$ and $W_{(j,\chi')}$ being adjacent if 
$i\not= j$ and $W_{(i,\chi)}\cap W_{(j,\chi')} = \emptyset$.
Note that $\name(G)$ is very close to the complement of the 
Eppstein's semicube graph Sc$(G)$. 

A path $P$ of $\name(G)$ with the property that 
$|P\cap \{W_{(i,0)}, W_{(i,1)}\}| \leq 1$ for each complementary 
pair of semicubes $W_{(i,0)}, W_{(i,1)}$, will be called a 
\emph{coordinating path}. 
A set of paths ${\cal P}$ of $\name(G)$ will be called a 
\emph{system of coordinating paths} provided that any 
$P\in {\cal P}$ is a coordinating path and for each 
complementary pair of semicubes $W_{(i,0)}, W_{(i,1)}$ there
is exactly one $P\in {\cal P}$ such that 
$|P\cap \{W_{(i,0)}, W_{(i,1)}\}| = 1$. 

\begin{lemma}
	\label{lem:upper}
	Let $G$ be a partial cube and let ${\cal P}$ be a system of coordinating 
	paths of $\name(G)$. Then there is an isometric embedding of $G$ into 
	$\Gamma_{f'}$, where 	$f'=\idim(G) + |{\cal P}| - 1\,.$
\end{lemma}

\proof
	Let $k=\idim(G)$, let $p = |{\cal P}|$, 
	and let ${\cal P} = \{P_1, \ldots, P_p\}$ be the given system of coordinating 
	paths of $\name(G)$. 	Let
	\begin{quote}
		$P_1: W_{(a_{1},\chi_{1})} \rightarrow W_{(a_{2},\chi_{2})} \rightarrow \cdots
		\rightarrow W_{(a_{i_1},\chi_{i_1})}$ \\
		$P_2: W_{(a_{i_1+1},\chi_{i_1+1})} \rightarrow W_{(a_{i_1+2},\chi_{i_1+2})} \rightarrow \cdots
		\rightarrow W_{(a_{i_2},\chi_{i_2})}$ \\
		\phantom{x}$\vdots$ \\
		$P_p: W_{(a_{i_{p-1}+1},\chi_{i_{p-1}+1})} \rightarrow W_{(a_{i_{p-1}+2},\chi_{i_{p-1}+2})} \rightarrow \cdots
		\rightarrow W_{(a_{i_p},\chi_{i_p})}$.
	\end{quote}
	As the paths meet exactly one of the complementary semicubes exactly once, 
	$i_p = k$. More precisely, there is a bijection $\phi:\{a_1,a_2,\ldots, a_{i_p}\} 
	\rightarrow [k]$ such that if $\phi(a_i) = j$
	then either $W_{(a_i,\chi_i)} = W_{(j,0)}$ or $W_{(a_i,\chi_i)} = W_{(j,1)}$ holds. 

	For any vertex $u$ of $G$ and any $i\in [k]$ set 
	\[
		\bar{u}^{(i)}=\begin{cases}
					1 & \mbox{ if $u\in W_{(a_i, \chi_i)}$;}\\
					0 & \mbox{ otherwise.}\\					
					\end{cases}
	\]
	Assigning the $k$-tuple 
	$$u=\bar{u}^{(1)}\bar{u}^{(2)}\ldots \bar{u}^{(i_p)}$$
	to any vertex $u$ of $G$ yields
	the canonical isometric embedding of $G$ into $Q_k = Q_{i_p}$. 
	Now assign to $u$ the following $d$-tuple: 
	$$\bar{u}^{(1)}\ldots \bar{u}^{(i_1)}0
	  \bar{u}^{(i_1+1)}\ldots \bar{u}^{(i_2)}0\ \ldots\ 0
	  \bar{u}^{(i_{p-1}+1)}\ldots \bar{u}^{(i_p)}\,.$$
	In this way, $G$ is embedded into $Q_{f'}$, where $f'=k+p-1$. Moreover, 
	the embedding is clearly still isometric. Because 
	$W_{(a_i,\chi_i)}\cap W_{(a_{i+1},\chi_{i+1})} = \emptyset$ provided 
	that $W_{(a_i,\chi_i)}$ and $W_{(a_{i+1},\chi_{i+1})}$ are connected
	by an edge of some path $P_j$, the labeling of $u$ 
	is a Fibonacci string. Hence we have described an isometric embedding
	of $G$ into $\Gamma_{f'}$.
\qed

Let $p(\name(G))$ be the minimum size of a system of coordinating paths 
of $\name(G)$. Then: 

\begin{theorem}
	\label{thm:glavni}
	Let $G$ be a partial cube. Then 
	$$\fdim(G) = \idim(G) + p(\name(G)) - 1\,.$$ 
\end{theorem}

\proof
	Let $p=p(\name(G))$, $k=\idim(G)$, and $f=\fdim(G)$. 
	If readily follows from Lemma~\ref{lem:upper}
	and the definition of $p(\name(G))$ that 
	$f \le k + p - 1$. 

	Consider now $G$ isometrically embedded into $\Gamma_f$. 
	For $u\in V(G)$ let 
	$u^{(1)}\ldots u^{(f)}$ be the embedded vertex. 
	Let $1\leq i_1 < i_2 < \cdots < i_r\le f$ be the indices 
	for which all the vertices of $G$ are labeled 0. 
	That is, $u^{(i_j)}=0$ holds for any $u\in V(G)$ and any 
	$i_j$,  $1\le j\le r$. Then 
	$$\beta(u) = u^{(1)}\ldots u^{(i_1-1)}u^{(i_1+1)}\ldots u^{(i_2-1)}
	u^{(i_2+1)}\ldots u^{(i_{r-1}-1)}u^{(i_{r-1}+1)}\ldots u^{(i_r)}$$
	is an isometric embedding into $Q_{f-r}$. 

	We next assert that for any coordinate $i$ of the $(f-r)$-tuples $\beta$, 
	$Y_i = \{\beta^{(i)}(u)\ |\ u\in V(G)\} = \{0,1\}$. 
	Note first that $Y_i\not= \{1\}$ because otherwise the $i$th coordinate 
	could be removed and hence we would isometrically embed $G$ into 
	$\Gamma_{f-1}$. On the other hand $Y_i\not= \{0\}$ since we have 
	removed all such coordinates in the construction of $\beta$. 
	Hence the assertion. 
	However, this implies that $f$ is an irredundant 
	embedding and therefore 
	$$k = \idim(G) = f - r\,.$$

	For a given coordinate $\ell$ of $\beta$, set 
	$W_\ell = \{u\in V(G)\ |\ \beta^{(\ell)}(u) = 1\}$. Then $W_\ell$ 
	is a semicube. Moreover, because $\beta$ is obtained from 
	Fibonacci strings, the paths
	\begin{quote}
	$W_{1} \rightarrow W_{2} \rightarrow \ldots \rightarrow W_{i_1-1}$, \\
	$W_{i_1+1} \rightarrow W_{i_1+2} \rightarrow \ldots \rightarrow W_{i_2-1}$, \\
	\phantom{x}$\vdots$ \\
	$W_{i_{r-1}+1} \rightarrow W_{i_{r-1}+2} \rightarrow \ldots \rightarrow W_{i_r}$,
	\end{quote}
	form a system of coordinating paths with $r+1$ paths. 
	Consequently, $r+1 \ge p$ and hence
	$$k = f - r \le f - p + 1\,.$$
	We conclude that $f\ge k + p - 1$ which completes the proof. 
\qed

Note that Proposition~\ref{prp:pc} also follows easily from 
Theorem~\ref{thm:glavni}.

\subsection{Particular cases}
\label{sec:combin-particular}

It is interesting to ask which partial cubes have extremal  
Fibonacci dimension. Interestingly, the minimum case is 
difficult; see Section~\ref{sec:bad}. However, there is a neat
characterization for the maximum case, which we provide next. 
Afterwards we establish the Fibonacci dimension of the Cartesian product
of graphs and the Fibonacci dimension of trees.

The \emph{crossing graph}
$G^\#$ of a partial cube $G$ has the $\Theta$-classes of $G$ as
its nodes, where two nodes of $G^\#$ are joined by an edge
whenever they cross as $\Theta$-classes in $G$; see~\cite{klmu-02}.
More precisely, if $W_{(a,0)}, W_{(a,1)}$ and $W_{(b,0)}, W_{(b,1)}$ are
pairs of complementary semicubes corresponding to 
$\Theta$-classes $E$ and $F$, then $E$ and $F$ cross if 
each semicube has a nonempty intersection with the semicubes from 
the other pair; that is, it holds that $W_{(a,0)}\cap W_{(b,0)}$, 
$W_{(a,0)}\cap W_{(b,1)}$, $W_{(a,1)}\cap W_{(b,0)}$, and $W_{(a,1)}\cap W_{(b,1)}$
are nonempty.

\begin{corollary}
\label{cor:extrem2}
	Let $G$ be a partial cube with  $\idim(G)=k$.  
	Then $\fdim(G)=2k-1$ if and only if $G^\# = K_k$. 
\end{corollary}

\proof
	By Theorem~\ref{thm:glavni}, $\fdim(G)=2k-1$ if and only
	if $p(\name(G)) = k$. This holds if and only if 
	$\name(G)$ has no edges which is in turn true if and 
	only if for any distinct $i,j\in [k]$ the semicubes
	$W_{(i,0)}$ and $W_{(i,1)}$ nontrivially intersect 
	$W_{(j,0)}$ and $W_{(j,1)}$. But this is true if and only if the corresponding 
	$\Theta$-classes cross. 
\qed

A characterization of complete crossing graphs in terms of 
the expansion procedure is given in~\cite{klmu-02}: 
$G^\#$ is complete if and only if $G$ can be obtained from 
$K_1$ by a sequence of all-color expansions. We also 
note that among median graphs only hypercubes have complete 
crossing graphs~\cite{mcmu-98}.

\begin{corollary}\label{co:cp}
	For any partial cubes $G$ and $H$, $\fdim(G\cp H) = \fdim(G)
	+ \fdim(H) + 1$.
\end{corollary}

\proof
	It is easy to infer that $\name(G\cp H)$ is isomorphic to $\name(G) \cup \name(H)$.
	Therefore, $p(\name(G\cp H)) = p(\name(G)) + p(\name(H))$.
	Since it is well-known that
	$\idim(G\cp H) = \idim(G) + \idim(H)$ we have
	\begin{eqnarray*}
	\fdim(G\cp H) & = & \idim(G\cp H) + p(\name(G\cp H)) - 1 \\
	& = & \idim(G) + \idim(H) + p(\name(G) \cup \name(H)) - 1 \\
	& = & \idim(G) + \idim(H) + p(\name(G)) + p(\name(H)) - 1 \\
	& = & (\idim(G) + p(\name(G)) - 1) + (\idim(H) + p(\name(H)) \\
	& = & \fdim(G) + \fdim(H) + 1\,,
	\end{eqnarray*}
	where for the first equality Theorem~\ref{thm:glavni} is applied.
\qed

\begin{corollary}\label{co:tree}
For any tree $T$, $\fdim(T) = \idim(T)=|E(T)|$.  
\end{corollary}

\proof
	Let $n=|V(T)|$. It is well-known that $\idim(T)=|E(T)|=n-1$,
	and that each edge $e$ of $T$ constitutes a $\Theta$-class~\cite{grwi-85} (cf.~\cite[Corollary 3.4.]{imkl-00}).
	This means that each edge $e\in E(T)$ defines a pair of complementary semicubes:
	each semicube is the set of vertices in one of the two subtrees of $T-e$.
		
\begin{figure}
	\centering
	\includegraphics[scale=1.2]{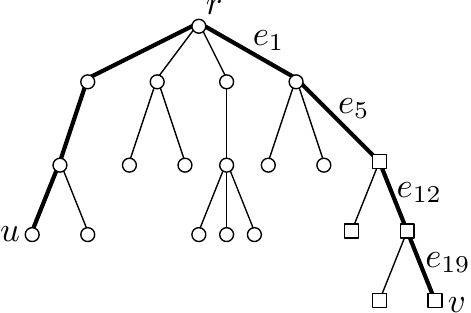}
	\caption{A tree with a longest path $P$ marked with thicker edges.
		In the proof of Corollary~\ref{co:tree}, $P_1$ would be the path between $r$ and $v$,
		$P_2$ would be the path between $r$ and $u$,
		the labeling of the edges of $P_1$ corresponds to a proper enumeration,
		and the nodes marked with squares correspond to the semicube $W_{e_5}$.}
	\label{fig:tree}
\end{figure}

	Let $P$ be a longest path in the tree $T$. We split $P$ at a vertex $r$
	into two subpaths $P_1,P_2$, 
	such that $P_1$ and $P_2$ have the same length (if $P$
	has an even number of edges), 
	or differ by one edge (if $P$ has an odd number of edges).
	Without loss of generality, let us assume 
	that $P_1$ is not strictly shorter than $P_2$.
	Therefore $|E(P_1)|=|E(P_2)|$ if $|E(P)|$ is even
	and $|E(P_1)|=1+|E(P_2)|$ if $|E(P)|$ is odd.
	See Figure~\ref{fig:tree}.
	
	It may be convenient to visualize $T$ as rooted at $r$.
	We further define the \emph{level} of an edge $xy$ of $T$ as the minimum
	of $d_T(r,x),d_T(r,y)$.
	For any edge $e$ of $T$, let $W_e$ denote the subset of vertices
	in the subtree $T-e$ that does not contain the vertex $r$. 
	As noted before, $W_e$ is a semicube, and hence a node of $\name(T)$, for any $e\in E(T)$.
	 
	Let $e_1,e_2,\dots e_{n-1}$ be an enumeration of the edges of $T$ 
	with the following properties:
	(a) any edge at level $i$ is listed before any edge at level $i+1$, and
	(b) the edge of $P_1$ at level $i$ is the first edge at level $i$ in the enumeration.
	Consider the sequence of semicubes $W_{e_1},W_{e_2},\ldots,W_{e_{n-1}}$.
	If the edges $e_i$ and $e_{i+1}$ are at the same level, 
	then clearly $W_{e_i}\cap W_{e_{i+1}}=\emptyset$. 
	If $e_i$ and $e_{i+1}$ are not at the same level, then $e_{i+1}$ must be an edge on $P_1$
	while $e_i$ cannot be an edge on $P_1$. Therefore we also have $W_{e_i}\cap W_{e_{i+1}}=\emptyset$
	in this case. This means that 
		$W_{e_1} \rightarrow W_{e_2} \rightarrow \ldots \rightarrow W_{e_{n-1}}$ 
	is a path in $\name(G)$, and furthermore forms a system of coordinating paths because
	it visits each complementary pair of semicubes exactly once.
	We conclude that $p(\name(T))=1$, and thus $\fdim(T)=\idim(T)$ by Theorem~\ref{thm:glavni}. 
\qed

\subsection{Relation to lattice dimension}

Using the lattice dimension $\ldim(G)$, we can provide
upper and lower bounds on the Fibonacci dimension $\fdim(G)$.
The first bound improves upon Proposition~\ref{prp:pc}.

\begin{figure}[t]
\centering\includegraphics[scale=0.5]{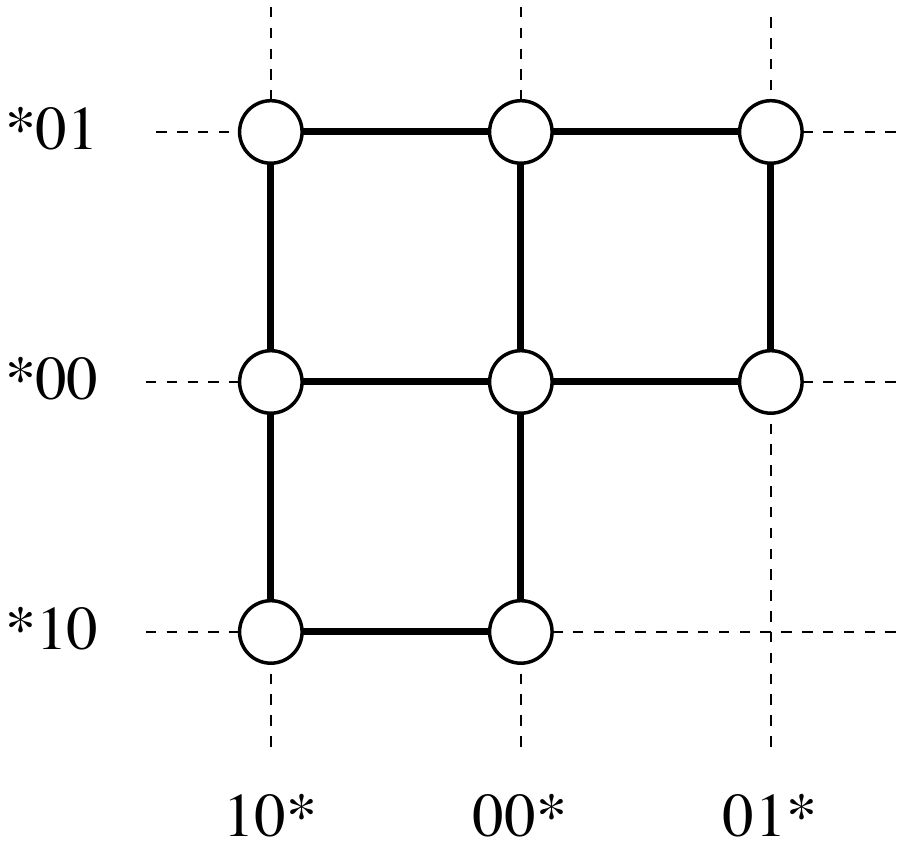}
\caption{A lattice embedding of $\Gamma_4$.}
\end{figure}

\begin{proposition}
\label{prp:ldim1}
	Let $G$ be a partial cube. Then 
		$\fdim\le \idim(G)+\ldim(G)-1$.
\end{proposition}
 
\proof
	For any integers $a,b$ with $a\le b$, 
	we use $P_{(a,b)}$ to denote the subgraph of $\mathbb{Z}^1$ 
	induced by vertices $a,a+1,\ldots, b-1,b$.
	Hence $P_{(a,b)}$ is a path on $b-a+1$ vertices and 
	$\fdim(P_{(a,b)}) = b-a$ by Corollary~\ref{co:tree}.
	
	Let $\ell=\ldim(G)$ and consider an isometric embedding $\beta$ of $G$ into $\mathbb{Z}^\ell$.
	For each coordinate $i\in [\ell]$, let $a_i=\min\{\beta^{(i)}(v)\mid v \in V(G)\}$
	and let $b_i=\max \{\beta^{(i)}(v)\mid v \in V(G)\}$.
	It is shown in~\cite[Lemma 1]{ep-05} that $\sum_i (b_i-a_i)$ is precisely $\idim(G)$.
	By the choice of $a_i,b_i$, the embedding $\beta$ 
	is also an isometric embedding of $G$ into the Cartesian product
	$P_{(a_1,b_1)}\cp P_{(a_2,b_2)}\cp \cdots \cp P_{(a_\ell,b_\ell)}$,
	and therefore
	\[
		\fdim(G)\le \fdim\left(P_{(a_1,b_1)}\cp P_{(a_2,b_2)}\cp \cdots \cp P_{(a_\ell,b_\ell)}\right).
	\]
	Since Corollary~\ref{co:cp} implies 
	\begin{align*}
		\fdim\left(P_{(a_1,b_1)}\cp P_{(a_2,b_2)}\cp \cdots \cp  P_{(a_\ell,b_\ell)}\right)
				\, &=\, \left(\sum_{i=1}^\ell \fdim(P_{(a_i,b_i)}) \right) + (\ell-1)\\
				&=\, \left(\sum_{i=1}^\ell (b_i-a_i) \right) + (\ell-1)\\
				&=\, \idim(G) + \ldim(G)-1,
	\end{align*}
	we conclude that $\fdim(G) \le \idim(G) + \ldim(G)-1$.
\qed

\begin{proposition}
\label{prp:ldim2}
	Let $G$ be a partial cube. Then 
	$\ldim(G)\le \lceil \fdim(G)/2\rceil$.
\end{proposition}

\proof
	Consider the Fibonacci cube $\Gamma_f$ for $f\ge 3$, and
	let $u^*$ denote the last $f-2$ entries of each tuple $u\in V(\Gamma_f)$.
	Define an embedding $\beta$ of $\Gamma_f$ into $\mathbb{Z}^1\cp \Gamma_{f-2}$ by 
	\[
		\beta(u)=\begin{cases}
					(0,u^*) & \mbox{ if $u=01u^*$;}\\
					(1,u^*) & \mbox{ if $u=00u^*$;}\\
					(2,u^*) & \mbox{ if $u=10u^*$.}\\					
					\end{cases}
	\]
	It is straightforward to see that $\beta$ is an isometric embedding.
	Using induction on the Fibonacci dimension, with base cases
	$\ldim(\Gamma_1)=\ldim(\Gamma_2)=1$, we obtain 
	\[
		\ldim(\Gamma_f)\,\le\, 1+\ldim(\Gamma_{f-2})\, \le\, 1+ \lceil (f-2)/2\rceil
					\,=\, \lceil f/2\rceil.
	\] 
	If a partial cube isometrically embeds into $\Gamma_f$, we
	then have $\ldim(G)\le \ldim(\Gamma_f)\le \lceil f/2\rceil$, and
	the result follows.
\qed

For graphs with low lattice dimension, we may determine the Fibonacci dimension exactly:

\begin{proposition}
Suppose that $\ldim(G)=2$. Then $\fdim(G)=\idim(G)+i$, where $i=1$ when $G$ is isomorphic to the Cartesian product of two paths and $i=0$ otherwise.
\end{proposition}

\proof
	When $G$ is isomorphic to the product of two paths, the result follows 
	from Corollaries~\ref{co:cp} and~\ref{co:tree}. 
	Otherwise, $G$ is a proper subgraph of $P_1\cp P_2$, where $P_1$ and $P_2$ are two paths with total length equal to $\idim(G)$. Among the four corner vertices of $P_1\cp P_2$ determined by pairs of endpoints of the two paths, at least one corner must be absent in $G$ if $G$ is to be a proper subgraph of the product of paths; we may assume without loss of generality that this missing corner corresponds to the last vertex of $P_1$ and the first vertex of $P_2$.
	
	We may embed $P_1$ isometrically into a Fibonacci cube (following Corollary~\ref{co:tree}) using the coordinates
$$101010\dots, \ 001010\dots, \ 000010\dots, \quad \dots, \quad \dots 010000, \ \dots 010100,\  \dots 010101$$
when $P_1$ has even length, or with a similar pattern when $P_1$ has odd length. That is, we start with an alternating sequence of zeros and ones, remove the ones one at a time, and then add ones one at a time to end with the opposite alternating sequence of ones and zeros. This pattern can be chosen in such a way that the final coordinate is zero for all vertices of $P_1$ except for its the last vertex. Similarly, we may embed $P_2$ isometrically into a set of Fibonacci strings in such a way that the initial coordinate is zero except in the first vertex of $P_2$. Concatenating these two representations of positions in $P_1$ and $P_2$ produces an irredundant isometric embedding of $G$ into a Fibonacci cube.
\qed

\section{Algorithmic aspects}
\label{sec:algor}

\subsection{Bad news}
\label{sec:bad}

We show that it is NP-complete to decide if the isometric and the Fibonacci dimension 
of a given graph is the same. 
Furthermore, we show that it is NP-hard to approximate the Fibonacci dimension within $(741/740)-\varepsilon$,
for any constant $\epsilon>0$.

Let $G$ be a graph with $n$ vertices. 
We assume for simplicity that $V(G)=[n]$,
and use $a,b$ to refer to the vertices of $G$.
Let $\bar G$ be the complementary graph of $G$.

\begin{lemma}
\label{lem:idimHG}
	Let $H$ be either the simplex graph $\simplex(G)$ or the 2-simplex graph $\simplex_2(G)$.
	Then $H$ is a partial cube with $\idim(H)=n$.
\end{lemma}

\begin{figure}
	\centering
	\includegraphics[scale=1]{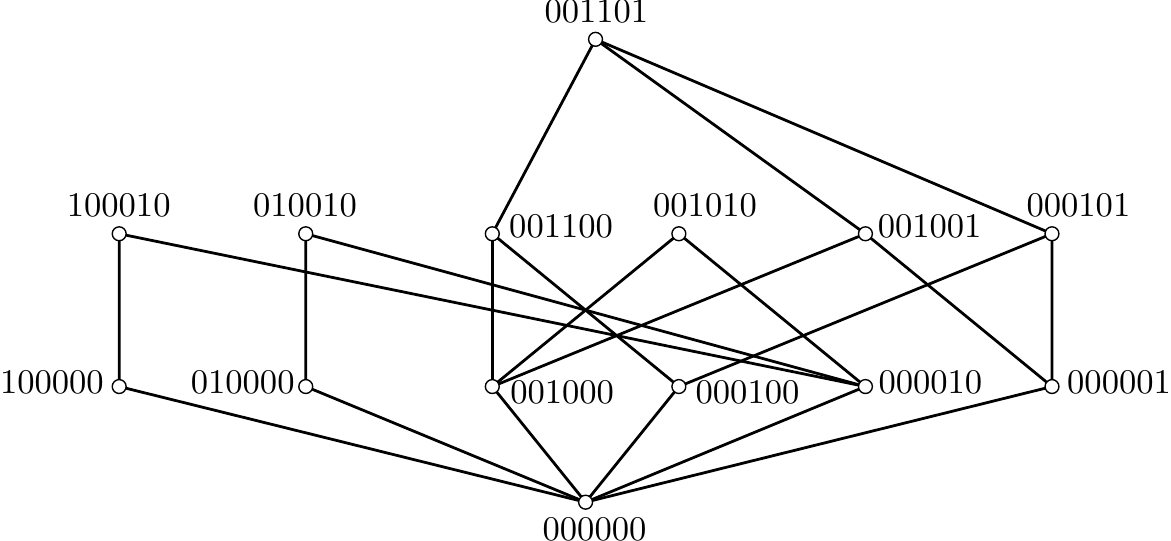}
	\caption{Isometric embedding of $\simplex(G)$ for the graph $G$ of Figure~\ref{fig:HG}, left.}
	\label{fig:isometric}
\end{figure}

\proof
	Consider the embedding 
	$\beta:H\rightarrow Q_n$ given as follows:
	\begin{itemize}
		\item for $u_\emptyset$ we set 
				$\beta(u_\emptyset)=u^{(1)}\ldots u^{(n)}$ with $u^{(i)}=0$ for all $i\in [n]$;
		\item for each $a\in [n]$ we set 
				$\beta(u_a)=u^{(1)}\ldots u^{(n)}$ with $u^{(a)}=1$ and $u^{(i)}=0$ for all $i\in [n]\setminus\{a \}$;
		\item for each node $u_K$ of $H$, we set 
				$\beta(u_{K})=\sum_{a\in K} \beta(u_a)$.
	\end{itemize}
	See Figure~\ref{fig:isometric} for an example when $H=\simplex(G)$.
	It is straightforward to see that $\beta$ is an isometric embedding of $H$ into $Q_n$, and
	hence $H$ is a partial cube. 
	Moreover, $\beta$ is irredundant:
	$\beta^{(i)}(u_a)$ is nonzero if and only if $a\not = i$.
	Since there is an irredundant isometric embedding of a graph $H$ into $Q_k$ 
	if and only if $\idim(H) = k$, it follows that $\idim(H)=n$.
\qed

In fact, stronger result that Lemma~\ref{lem:idimHG} was proved in~\cite{Banvan-PotLMS-89} for 
$\simplex(G)$ and in~\cite{imklmu-99} for $\simplex_2(G)$: $H$ is a median graph. 

\begin{lemma}
\label{lem:Xsimplex}
	Let $H$ be either the simplex graph $\simplex(G)$ or the 2-simplex graph $\simplex_2(G)$.
	There is a set $\W$ of semicubes of $H$ with the following properties:
	\begin{itemize}
		\item[\textup{(a)}] Each node in $\W$ has degree zero in $\name(H)$.
		\item[\textup{(b)}] Each pair of complementary semicubes of $H$ has a node
			in $\W$.
		\item[\textup{(c)}] $\name(H)-\W$ is isomorphic to $\bar G$.
	\end{itemize}
\end{lemma}

\begin{figure}
	\centering
	\includegraphics[width=.9\textwidth]{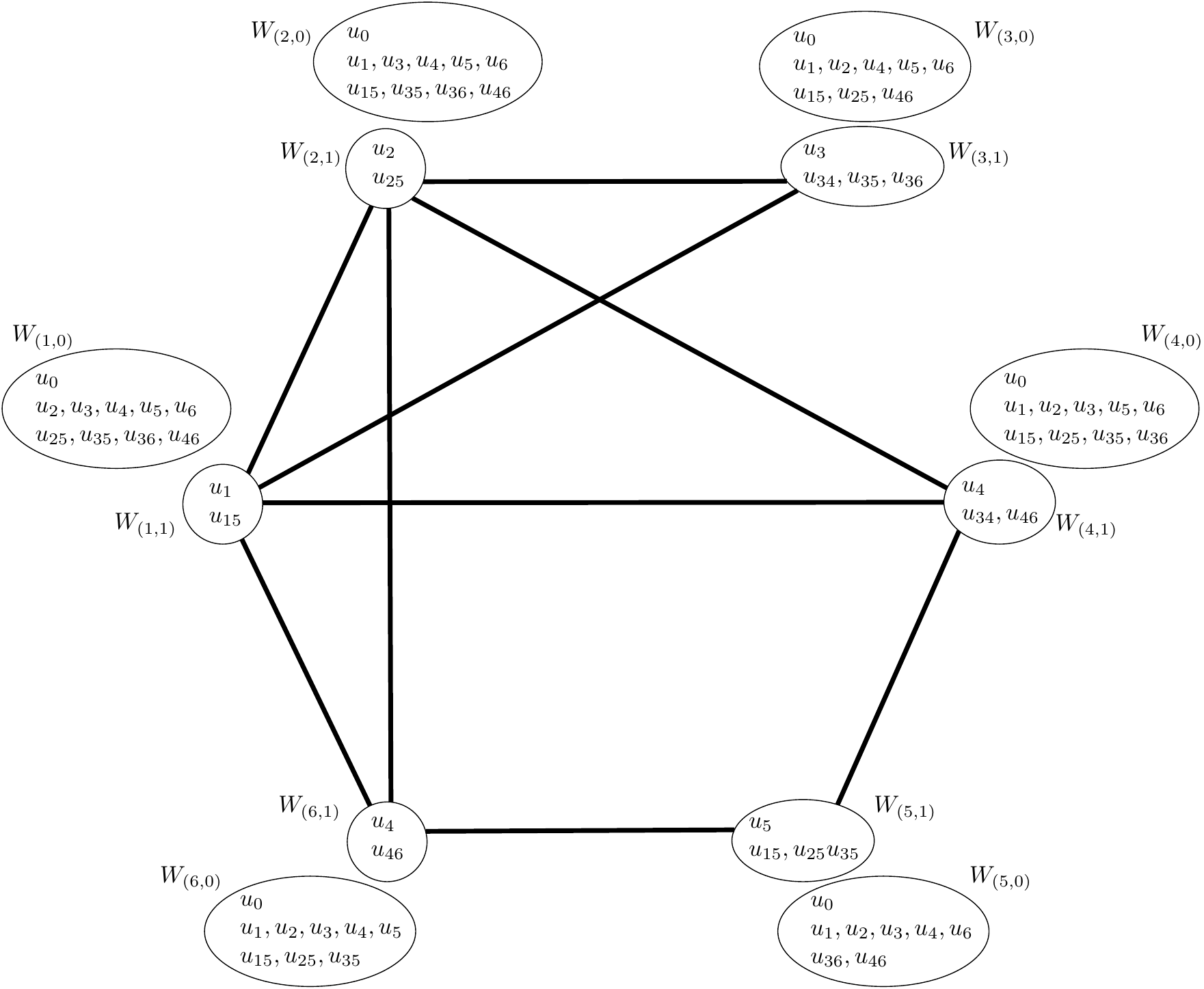}
	\caption{The graph $\name(\simplex_2(G))$ for the graph $G$ of Figure~\ref{fig:HG}, left.}
	\label{fig:XHG}
\end{figure}

\proof
	We will use the isometric embedding $\beta$ given in the proof of Lemma~\ref{lem:idimHG}.
	For any $a\in [n]$ we then have the semicubes
	\begin{align*}
		W_{(a,0)}\, 
				&= \, \{ u_K\in V(H) \mid \beta^{(a)}(u_K)=0 \} \\
				&= \, \{ u_K\in V(H) \mid \mbox{$a$ is not a vertex in $K$} \},
	\end{align*}
	and
	\begin{align*}
		W_{(a,1)}\, 
				&= \, \{ u_K\in V(H) \mid \beta^{(a)}(u_K)=1\} \\
				&= \, \{ u_K\in V(H) \mid \mbox{$a$ is a vertex in $K$} \}.
	\end{align*}
	Let us now consider the graph $\name(H)$. See Figure~\ref{fig:XHG} for an example.
	The node set of $\name(H)$ is
	$W_{(a,\chi)}$, $(a,\chi)\in [n]\times \{0,1\}$. For the edge set, we have the following
	properties:
	\begin{itemize}
		\item There is no edge between $W_{(a,0)}$ and $W_{(b,0)}$ because $u_\emptyset \in W_{(a,0)} \cap W_{(b,0)}$.
		\item There is no edge between $W_{(a,0)}$ and $W_{(b,1)}$ because $u_b \in W_{(a,0)}\cap W_{(b,1)}$.
		\item There is an edge between $W_{(a,1)}$ and $W_{(b,1)}$ if and only if $ab\notin E(G)$. Indeed,
			there is a vertex $u_K$ of $H$ in $W_{(a,1)} \cap W_{(b,1)}$ if and only if $a$ and $b$
			are vertices in the clique $K$, which happens precisely when $ab$ is an edge of $G$.
			Therefore $W_{(a,1)} \cap W_{(b,1)} \not = \emptyset$
			if and only if $ab\in E(G)$.
	\end{itemize}
	It follows that each node $W_{(a,0)}$, $a\in [n]$, has degree zero in $\name(H)$.
	Therefore, the subfamily of nodes $\W= \{ W_{(a,0)}\mid a\in [n]\}$ of $\name(H)$
	satisfies properties (a) and (b) in the lemma. The graph $\name(G)-\W$ contains
	only the nodes $W_{(a,1)}$, $a\in [n]$. 
	Since there is an edge between $W_{(a,1)}$ and $W_{(b,1)}$ if and only if $ab\notin E(G)$,
	the mapping $a\mapsto W_{(a,1)}$ is an isomorphism between $\bar G$, and property (c) follows.
\qed

Let $(1,2)$-TSP denote the (metric) symmetric Traveling Salesman Problem 
in which all distances are either 1 or 2.
The $(1,2)$-TSP problem is NP-hard. Furthermore, 
Engebretsen and Karpinski~\cite{enka-06} have shown that
it is NP-hard to approximate the $(1,2)$-TSP within $(741/740)-\varepsilon$
for every constant $\varepsilon>0$.
On the positive side,
Berman and Karpinski~\cite{beka-06} have given an $(8/7)$-approximation algorithm for $(1,2)$-TSP.

Any graph $G$ naturally defines an instance $I_G$ of $(1,2)$-TSP, where the points of the metric
space are the vertices of $G$, and the distance between two points is 
1 if there is an edge between them in $G$, and 2 otherwise.
Let $\ell(I_G)$ denote the length of the optimal tour for an instance $I_G$ of $(1,2)$-TSP.

For later use, it will be convenient to exchange now the roles of $G$ and its complementary graph $\bar G$.
\begin{lemma}
\label{lem:tsp}
	Let $H$ be either the simplex graph $\simplex(\bar G)$ or the 2-simplex graph $\simplex_2(\bar G)$.
	The graph $G$ has a Hamiltonian path if and only if $\fdim(H)=n$.
	If $G$ does not have a Hamiltonian path, then $\fdim(H)=\ell(I_{G})-1$.
\end{lemma}

\proof
	Consider the set of nodes $\W$ in $\name(H)$ given by Lemma~\ref{lem:Xsimplex}.
	Since each node of $\W$ has degree zero in $\name(H)$ and $\W$ contains one
	semicube from each pair of complementary semicubes of $H$, we can
	just disregard the nodes $\W$ for finding the value $p(\name(H))$.
	When we disregard the nodes $\W$, we obtain $\name(H)-\W$, which is isomorphic
	to $G$ because of property (c) in Lemma~\ref{lem:Xsimplex}.
	(Recall we exchanged the roles of $G$ and $\bar G$.)
	It follows that $p(\name(H))$ is the minimum number of vertex-disjoint paths that 
	are needed to cover each vertex of $G$.
	
	The graph $G$ has a Hamiltonian path if and only if $p(\name(H))=1$. 
	Using Theorem~\ref{thm:glavni} and Lemma~\ref{lem:idimHG},
	this is equivalent to 
	\[
		\fdim(H)\,=\,\idim(H)+ p(\name(H))-1 \,=\, \idim(H)=n.
	\]
	
	If $G$ does not have a Hamiltonian path, then
	$\ell(I_{G})$ is $|V(G)|=n$ 
	plus the minimum number of vertex-disjoint
	paths that are needed to cover each vertex of $G$.
	Thus $\ell(I_{G})=n+p(\name(H))$.
	Using Theorem~\ref{thm:glavni} and Lemma~\ref{lem:idimHG} we conclude that 
	\[
		\fdim(H)\,=\,\idim(H)+ p(\name(H))-1 \,=\, n +p(\name(H))-1\,=\, \ell(I_{G})-1.\]
\qed

We next show that computing the Fibonacci dimension, or even to approximate it, is NP-hard. 

\begin{theorem}
\label{thm:hard1}
	It is NP-complete to decide if $\idim(H)=\fdim(H)$ for a given graph $H$.
\end{theorem}

\proof 
	Note that $\idim(H)$ can be computed in polynomial time~\cite{ah-95,ep-08,imkl-93}.
	Therefore, an explicit isometric embedding of $H$ into $\Gamma_{\idim(H)}$ would 
	be enough to check in polynomial time that $\idim(H)=\fdim(H)$. It follows that 
	the problem is in the class NP.
	
	To show hardness, consider the graph $H=\simplex_2(\bar G)$.
	It is clear that $H$ can be constructed in polynomial time for any given graph $G$;
	this is not necessarily true for $\simplex(\bar G)$ if $\bar G$ has large cliques.
	Lemma~\ref{lem:tsp} 
	implies that $\idim(H)=\fdim(H)$ 
	if and only if $G$ has a Hamiltonian path.
	Since deciding whether a graph has a Hamiltonian path is NP-complete~\cite{GJ-book}, 
	it is NP-hard to decide whether $\idim(H)=\fdim(H)$.
\qed

\begin{theorem}
\label{thm:hard2}
	It is NP-hard to approximate the Fibonacci dimension of a graph within $(741/740)-\varepsilon$
	for every constant $\varepsilon>0$.
\end{theorem}

\proof
	Assume that there is a constant $\varepsilon>0$
	and a polynomial time algorithm {\sc ApproxFib} that, 
	for any input graph $H$, computes a value $f'(H)$ such that 
	\[
		\fdim(H)\,\le\, f'(H)\,\le\, \left( \tfrac{741}{740}-\varepsilon\right)\fdim(H).
	\]
	Given any graph $G$ with $n$ vertices,
	we can apply algorithm {\sc ApproxFib} to the graph $H=\simplex_2(\bar G)$ to obtain
	a value $f'$ that satisfies
	\begin{equation}
		\label{eq:f'}
		\fdim(H)\,\le\, f'\,\le\, 
		\left( \tfrac{741}{740}-\varepsilon\right)\fdim(H).
	\end{equation}
	Consider the value $\ell'= f'+1$ as an approximation to $\ell(I_{G})$.
	
	From Lemma~\ref{lem:tsp} it follows that 
	\begin{equation}
		\label{eq:f}
		\ell(I_{G})-1 \, \le \, \fdim(H) \,\le\, \ell(I_{G}).
	\end{equation}
	(There is the special case when $G$ has a Hamiltonian cycle because then
	$\ell(I_{G})=n=\fdim(H)$.)
	Combining inequalities (\ref{eq:f'}) and (\ref{eq:f}) we obtain
	\[
		\ell(I_{G})\,\le \, \fdim(H)+1 \,\le\, f' +1 \,= \, \ell',
	\]
	and 
	\begin{align*}
		\ell' \,&= \, f'+1 \\
				& \le\, \left( \tfrac{741}{740}-\varepsilon\right)\fdim(H) +1\\
				& \le\, \left( \tfrac{741}{740}-\varepsilon\right)\ell(I_{G}) +1\\
				&=\, \left( \tfrac{741}{740}+\tfrac{1}{\ell(I_{G})}-\varepsilon\right)\ell(I_{G})\\
				&\le\, \left( \tfrac{741}{740}+\tfrac{1}{n}-\varepsilon\right)\ell(I_{G}).
	\end{align*}
	Since $2/\varepsilon$ is a constant, we may assume that $G$ has more than $2/\varepsilon$ vertices,  
	and thus
	\[
		\ell(I_{G}) \,\le\, \ell' \,\le \, \left( \tfrac{741}{740}+\tfrac{1}{n}-\varepsilon\right)\ell(I_{G})
			\,\le\, \left( \tfrac{741}{740}-\tfrac{\varepsilon}{2}\right)\ell(I_{G}).
	\]
	We then conclude that $\ell'$ can be computed in polynomial time
	and approximates the value $\ell(I_{G})$ within $(741/740)-(\varepsilon/2)$.
	However, Engebretsen and Karpinski~\cite{enka-06} have shown 
	that it is NP-hard to approximate the $(1,2)$-TSP within $(741/740)-\delta$
	for every constant $\delta>0$.
	Therefore, it is also NP-hard to approximate the 
	Fibonacci dimension of a graph within $(741/740)-\varepsilon$
	for every constant $\varepsilon>0$.
\qed
 
\subsection{Good news}
\label{sec:good}

We first provide an exact algorithm to compute $\fdim(G)$
whose running time is exponential in $\idim(G)$.
We then provide a (3/2)-approximation algorithm for arbitrary graphs,
and better approximation algorithms specialized to simplex graphs.

We assume that our input is a partial cube $G$ with $n$ vertices and also
that we are given an embedding $\beta$ of $G$ into $Q_{k}$, where $k=\idim(G)$. 
Such embedding can be constructed in $O(n^2)$ time~\cite{ep-08}\footnote{%
The algorithm in~\cite{ep-08} assumes the word-RAM model of computation. 
Without bit-manipulation, there
are algorithms~\cite{ah-95,imkl-93} taking $O(n^2\log n)$ time.}.
We first describe how to construct $\name(G)$ and then give
an algorithmic counterpart of Lemma~\ref{lem:upper}.

\begin{lemma}
\label{lem:XG}
	The graph $\name(G)$ can be computed in $O(k^2 n)$ time.
\end{lemma}
\proof
	Each semicube $W_{(i,\chi)}$ is identified by a pair $(i,\chi)\in [k]\times \{0,1\}$.
	We first construct the complete graph on the node set $\{ (i,\chi)\in [k]\times \{0,1\}\}$
	and then, for each vertex $v\in V(G)$, the edges $(i,\beta^{(i)}(v))(j,\beta^{(j)}(v))$ are removed for
	all distinct $i,j\in [k]$. The resulting graph is (isomorphic to) $\name(G)$.
	Using any standard data structure for graphs, 
	each edge can be deleted in constant time.
	For each of the $n$ vertices of $G$, we thus spend $O(k^2)$ time,
	for a total of $O(k^2 n)$ time.
\qed

\begin{lemma}
	\label{lem:upper2}
	Assume we are given a system of $p$ coordinating paths of $\name(G)$. 
	Then we can compute in $O(kn)$ time 
	an isometric embedding of $G$ into $\Gamma_{f'}$, where
	$f'=k + p - 1\,.$
\end{lemma}

\proof
	The proof given in Lemma~\ref{lem:upper} is constructive and can be implemented in $O(n(k+p))=O(kn)$ time.
\qed

From $\name(G)$ it is possible to compute $p(\name(G))$, and thus $\fdim(G)$,
in roughly $O(k!)$ time by trying all permutations of the indices $[k]$ to
obtain systems of coordinating paths of $\name(G)$.
We next improve this to a dependency that is exponential in $k$.

\begin{proposition}
\label{prp:exact}
	Given a partial cube $G$ with $n$ vertices and an isometric embedding $G\rightarrow Q_{k}$, 
	where $k=\idim(G)$,
	we can compute in $O(2^k k^2+ k^2 n)$ time an isometric embedding of $G$ 
	in $\Gamma_{f}$, where $f=\fdim(G)$.
\end{proposition}

\proof
	Firstly, we construct the graph $\name(G)$ using Lemma~\ref{lem:XG} in $O(k^2n)$ time.
	Secondly, we find in $O(2^k k^2)$ time
	a system of coordinating paths of $\name(G)$ with minimum size
	using dynamic programming, as described below.
	Finally, we use Lemma~\ref{lem:upper2} to construct the embedding in $O(kn)$ time.
	We only have to describe the second step.
	
	We compute the value $p(\name(G))$ using dynamic programming across subsets of 
	pairs of complementary semicubes, as follows. 
	Our approach is essentially the same as a standard one for TSP~\cite{be-60,heka-62}.
	For any subset of indices $I\subseteq [k]$,
	let $\name_I(G)$ denote the subgraph of $\name(G)$ induced by
	nodes $W_{(i,0)}, W_{(i,1)}$, $i\in I$.
	For any triple $(I,j,\chi)\in 2^{[n]}\times [n]\times \{0,1\}$ with $j\in I$,
	let $\pi(I,j,\chi)$ denote the minimum number of paths in
	a system of coordinating paths for $\name_I(G)$, with the property
	that $W_{(j,\chi)}$ is an end-node of some coordinating path.
	That is, $\pi(I,j,\chi)$ is the minimum
	number of paths in $\name(G)$ that visit each pair of complementary semicubes 
	$W_{(i,0)}, W_{(i,0)}$, $i\in I$, exactly once,
	it has one path ending at node $W_{(j,\chi)}$, and does not visit any
	semicube $W_{(i,*)}$ for $i\notin I$.
	
	It is clear that for any $j\in [k]$ it holds $\pi(\{j\},j,0)=\pi(\{j\},j,1)=1$.
	For subsets $I$ with more than one index there are two cases
	to distinguish, depending
	on whether the paths defining $\pi(I,j,\chi)$ have $W_{(j,\chi)}$
	as an isolated node or not. Therefore it holds
	\[
		\pi(I,j,\chi) = \min \begin{cases}
								\displaystyle 
								1 + \min_{(j',\chi')\in (I\setminus\{ j \})\times \{0,1\}}
								\;\pi(I\setminus \{j \},j', \chi')\\
								\quad \\
								\displaystyle \min_{\begin{minipage}{4.7cm}\centering \small
								$(j',\chi')\in (I\setminus\{ j \})\times \{0,1\}$ s.t.\\
								$W_{(j,\chi)}W_{(j',\chi')} \in E(\name(G))$
								\end{minipage}}
								\;\pi(I\setminus \{j \},j', \chi')
							\end{cases}
	\]
	Note that $\pi(I,j,\chi)$ only depends on values $\pi(I',j',\chi')$ with $|I'|=|I|-1$.
	Therefore, we can compute all values $\pi(I,j,\chi)$ by considering them for increasing values
	of $|I|$, at a cost of $O(k)$ per value. Since there are at most $|2^{[k]}|\cdot |[k]|\cdot 2=O(2^k k)$
	tuples $(I,j,\chi)$ to consider, we can compute in $O(2^k k^2)$ time the values $\pi(I,j,\chi)$
	for all tuples $(I,j,\chi)\in 2^{[n]}\times [k]\times \{0,1\}$ with $j\in I$.
	Finally, it holds that 
	\[
		p(\name(G))= \min_{(j,\chi)\in [n]\times \{0,1\} } \pi([n],j,\chi),
	\]
	and hence we can recover $p(\name(G))$ in $O(k)$ time. To obtain the actual system
	of coordinating paths, we only need to augment each entry $(I,j,\chi)$ with
	a list of the paths that define $\pi(I,j,\chi)$.
\qed

We next move onto approximation algorithms.
Note that the value
$\idim(G)+\ldim(G)-1$ is a $(3/2)$-approximation to the value $\fdim(G)$
because of Propositions~\ref{prp:ldim1} and~\ref{prp:ldim2}.
Moreover, the proofs of Propositions~\ref{prp:ldim1} and~\ref{prp:ldim2}
are constructive, and therefore we can use isometric embeddings of $G$
into $Q_{\idim(G)}$ and into $\mathbb{Z}^{\ldim(G)}$ to construct
an isometric embedding of $G$ into $\Gamma_{\idim(G)+\ldim(G)-1}$.
Since an isometric embedding of $G$ into $\mathbb{Z}^{\ldim(G)}$ can be computed
in polynomial time~\cite{ep-06}, we can then compute in polynomial time an embedding
of $G$ into $\Gamma_{f'}$ for $f'\le (3/2)\fdim(G)$.
We next give an alternative algorithm with the same performance (time and approximation factor)
that does not make the detour through 
finding an isometric embedding into $\mathbb{Z}^{\ldim(G)}$.

\begin{theorem}
	\label{thm:approx}
	Given a partial cube $G$ with $n$ vertices and an isometric embedding $G\rightarrow Q_{k}$, 
	where $k=\idim(G)$,
	we can compute in $O(k^2 n)$ time an isometric embedding of $G$ 
	in $\Gamma_{f'}$, where $f'\le (3/2)\fdim(G)$.
\end{theorem}

\proof
	We first describe the algorithm, then derive its running time,
	and finally discuss the bound on the dimension of the computed embedding.
	
	The algorithm is as follows.
	Firstly, construct the graph $\name(G)$.
	Secondly, construct the graph $Y(G)$ obtained from $\name(G)$ by identifying 
	each pair of complementary semicubes into a single node.
	Hence, $Y(G)$ has $k$ nodes.
	Thirdly, construct a matching $M_Y$ in $Y(G)$ of maximum cardinality.
	Let $M_\name$ denote the matching in $\name(G)$ that corresponds to $M_Y$.
	We can then regard each edge of $M_\name$ as a path in $\name(G)$ that
	passes through two nodes. Let $P_1,\ldots, P_{|M_\name|}$ denote these paths. 
	There are precisely $k-2|M_\name|$ pairs of complementary semicubes that are
	not adjacent to $M_\name$. For each of those pairs, we make a path consisting
	of a single semicube of the pair. This gives a family of 
	$|M_\name|+(k-2|M_\name|)= k- |M_\name|$ paths that form a 
	system of coordinating paths of $\name(G)$.	
	Finally, we compute the embedding into $\Gamma_{f'}$ given by Lemma~\ref{lem:upper2},
	where 
	\begin{equation}\label{eq:d'}
		f'= k+ (k- |M_\name|) -1= 2k - |M_\name| -1.
	\end{equation}
	This finishes the description of the algorithm. 
	Clearly, this algorithm computes a valid embedding of $G$ into $\Gamma_{f'}$.
	
	To derive its running time, note that $\name(G)$ is constructed in $O(k^2 n)$ time by
	Lemma~\ref{lem:XG}.
	We can then construct $Y(G)$ by identifying the nodes $(i,0)$ and $(i,1)$ of $\name(G)$
	for all $i\in [d]$.
	Finding a maximum matching $M_Y$ in $Y(G)$ takes $O(k^{5/2})$ time using~\cite{MicVaz} 
	because $Y(G)$ has $k$ nodes and $O(k^2)$ edges. From $M_Y$ we can recover the matching $M_\name$
	in $\name(G)$, and construct the embedding using Lemma~\ref{lem:upper2} in $O(kn)$ time.
	We conclude that the algorithms takes $O(k^{5/2}+ k^2 n)$ time, which is $O(k^2n)$
	because $k\le n$.
	
	It remains to bound $f'$. Let $p=p(\name(G))$ and let $P_1,P_2,\ldots P_p$ 
	be a system of coordinating paths of $\name(G)$.
	Consider these paths in $Y(G)$, 
	and let $E_P$ denote the set of edges appearing in $P_1,P_2,\ldots P_p$. It holds
	that 
	\begin{equation}\label{eq:p}
		p=k-|E_P|.
	\end{equation}
	Taking each other edge in the path $P_i$, we see that $P_i$
	contains a matching	with $\lceil |E(P_i)|/2 \rceil$ edges.
	Thus the paths $P_1,\ldots P_p$ in $Y(G)$ contain a matching with at least $|E_P|/2$ edges.
	We conclude that
	\begin{equation}\label{eq:M}
		|M_\name| \, = \, |M_Y| \,\ge \, \frac{|E_P|}{2}.
	\end{equation}
	Combining equations (\ref{eq:d'})-(\ref{eq:M}) we obtain
	\begin{align*}
		f' \, &= \, 2k - |M_\name| -1 \\
		  \,&\leq \, 2k - \frac{|E_P|}{2}  -1 \\
		  \,&= \, 2k - \frac{k-p}{2} -1 \\
		  \,&= \, \frac{3}{2} (k+p-1) - p +\frac{1}{2} \\
		  \,&\leq \, \frac{3}{2} (\fdim(G) ),
	\end{align*}
	where in the last step we have used Theorem~\ref{thm:glavni}.
\qed

We now turn our attention to approximation algorithms for simplex graphs.

\begin{theorem}
	\label{thm:approxsimplex}
	Let $\varepsilon \in (0,1)$ be a constant.
	Given a simplex graph $H$ with $n$ vertices and an isometric embedding $H\rightarrow Q_{k}$, 
	where $k=\idim(H)$,
	we can compute in $O(n^{2+ 1/\varepsilon})$ time an isometric embedding of $H$ 
	in $\Gamma_{f'}$, where $f'\le (1+\varepsilon)\fdim(H)$.
\end{theorem}

\proof
	Since $H$ is a simplex graph, then $H=\simplex(G)$ for some graph $G$.
	It holds $k=\idim(H)=|V(G)|$ because of Lemma~\ref{lem:idimHG}.
	The graph $\name(H)$ can be constructed in $O(k^2 n)$ time by Lemma~\ref{lem:XG}.	
	We can construct the set of nodes $\W$ of Lemma~\ref{lem:Xsimplex} by placing in $\W$,
	for each $i\in [k]$, either $W_{(i,0)}$ or $W_{(i,1)}$, whichever has degree zero in $\name(G)$.
	By property (c) of Lemma~\ref{lem:Xsimplex} it holds that $\name(H)-\W$ and $\bar G$ are isomorphic.
	
	We construct a system of coordinating paths of $\name(H)$ with the following greedy procedure.
	We start with a system of $k$ coordinating paths $\mathcal{P}=\{ P_1,\ldots ,P_k\}$, 
	where each path consists of a node from $V(\name(H))\setminus \W$. 
	Then, we repeat the following step as many times as possible:
	if there are paths $P,P'\in \mathcal{P}$ that can be joined by adding
	an edge between two of its extreme nodes, we do so, and replace in $\mathcal{P}$ the
	paths $P,P'$ by the new path. 
	This step may be repeated at most $k-1$ times, in which case we end up with a single path in $\mathcal{P}$. 
	Since each repetition of the step takes $O(k^2)$ time,
	the whole procedure needs $O(k^3)$ time.
	Let $\mathcal{P}$ denote the resulting system of coordinating paths.	

	We distinguish two cases depending on the size of $\mathcal{P}$.
	If $|\mathcal{P}|\le \varepsilon k$, then we can use Lemma~\ref{lem:upper2}
	with $\mathcal{P}$ to construct an embedding of $H$ into $\Gamma_{f'}$,
	where $f'\le k+|\mathcal{P}| -1 \le (1+\varepsilon) k \le (1+\varepsilon) \fdim(H)$.
	If $|\mathcal{P}|> \varepsilon k$, then by selecting an extreme node in each path
	of $\mathcal{P}$ we obtain in $\name(G)-\W$ an independent set of nodes of cardinality
	at least $\varepsilon k$. Since $\name(G)-\W$ and $\bar G$ are isomorphic,
	this means that $\bar G$ has an independent set with at least $\varepsilon k$ vertices,
	and hence $G$ has a clique $K$ with at least $\varepsilon k$ vertices.
	Since the input graph $H$ is the simplex graph of $G$, we conclude
	that $H$ has at least $2^{|K|}\ge 2^{\varepsilon k}$ vertices;
	that is $n\ge 2^{\varepsilon k}$.
	Using Proposition~\ref{prp:exact}, we can then compute an isometric
	embedding of $H$ into $\Gamma_{\fdim(H)}$ in 	
	\[
		O(2^k k^2+ k^2 n)= O\left( \left(2^{\varepsilon k}\right)^{1/\varepsilon} n^2+ n^3\right)= 
		O\left( n^{1/\varepsilon} n^2+ n^3\right)=O\left( n^{2+1/\varepsilon}\right)
	\]
	time. In either case, we obtain in $O\left( n^{2+1/\varepsilon}\right)$ time
	an isometric embedding of $H$ into $\Gamma_{f'}$, where $f'\le (1+\varepsilon)\fdim(H)$.	
\qed

Note that the last result takes time polynomial in the size of the given simplex graph.
However, we could consider that the simplex graph $\simplex(G)$ is described by $G$,
and give $G$ as the input.
The following result approximates the Fibonacci dimension of $\simplex(G)$ for a given graph $G$
in time that is polynomial in the size of $G$. 
Note that we cannot compute an explicit isometric
embedding of $\simplex(G)$ in polynomial time of because the size of $\simplex(G)$ may be exponential in $G$.

\begin{theorem}
	\label{thm:approxsimplex2}
	Let $\varepsilon \in (0,1)$ be a constant.
	Given a graph $G$ with $n$ vertices, 
	we can compute in polynomial time a value $f'$ 
	such that $\fdim(\simplex(G))\le f'\le (\tfrac{8}{7}+\varepsilon)\fdim(\simplex(G))$.
\end{theorem}
\proof
	If $G$ has less than $1/\varepsilon$ vertices, which is a constant,
	we can then compute $\fdim(\simplex(G))$ exactly in constant time.
	Let us assume henceforth that $G$ has at least $1/\varepsilon$ vertices.
	
	From Lemma~\ref{lem:tsp} it follows that 
	$\fdim(\simplex(G)) = \ell(I_{\bar G})-1$ if $\bar G$ does not have a Hamiltonian cycle,
	and $\fdim(\simplex(G)) = n = \ell(I_{\bar G})$ if $\bar G$ has a Hamiltonian cycle.
	Berman and Karpinski~\cite{beka-06} describe a polynomial-time
	$(8/7)$-approximation algorithm for $(1,2)$-TSP.
	Let $\ell'$ be the $(8/7)$-approximation to the value $\ell(I_{\bar G})$ returned by
	their algorithm.
	If $\ell'=n$, then $\bar G$ has a Hamiltonian cycle, and $\fdim(\simplex(G)) = n$.
	If $\ell'> n$, consider $f'=\ell'-1$ as an approximation to $\fdim(\simplex(G))$.
	For $f'$ we have the upper bound 
	\[
		f' \,=\, \ell'-1 \, \le\, \tfrac{8}{7} \ell(I_{\bar G}) -1\,\le\,
			\tfrac{8}{7}(\ell(I_{\bar G})-1) + \tfrac{1}{7} \,\le\,
			\tfrac{8}{7}\fdim(\simplex(G)) + \tfrac{1}{7},	
	\]
	which using that $n\varepsilon \ge 1$ and $n\le \fdim(\simplex(G))$ leads to
	\[
		f' \,\le\, \tfrac{8}{7}\fdim(\simplex(G)) + \tfrac{\varepsilon n}{7} \,\le\,
			\tfrac{8}{7}\fdim(\simplex(G)) + \varepsilon \fdim(\simplex(G)) \,\le\,
			(\tfrac{8}{7}+\varepsilon)\fdim(\simplex(G)).	
	\]
	On the other hand, for $f'$ we also have the lower bound
	\[
		\fdim(\simplex(G))\,=\, \max\{n,\ell(I_{\bar G})-1\} \,\le\, \max\{n,\,\ell'-1 \}
		\, = \, \ell'-1 = f'.
	\]
	The result follows.
\qed

By Lemma~\ref{lem:tsp} it holds that $\fdim(\simplex(G))=\fdim(\simplex_2(G))$, 
and hence the same result holds
for the 2-simplex graph $\simplex_2(G)$.




\section*{Acknowledgments}

Work of D. Eppstein was supported in part by NSF grant 0830403 and by the Office of 
Naval Research under grant N00014-08-1-1015. Work of S. Cabello and S. Klav\v zar 
was supported in part by the Slovenian Research Agency, program P1-0297.
We also would like to thank David Johnson and Thore Husfeldt for clarifying the history 
of using dynamic programming for TSP. 


\bibliographystyle{abbrv}
\bibliography{biblio-fibonacci}

\end{document}